\documentclass[12pt,draftcls,onecolumn]{IEEEtran}
\usepackage{amsmath} 
\usepackage{amsfonts,mathrsfs}
\usepackage{amssymb}
\usepackage{color}
\usepackage{graphicx}
\usepackage{epsfig}
\usepackage{subfigure}

\newtheorem{lemma}{Lemma}{\bf}{}

\newtheorem{assumption}{Assumption}{\bf}{}

\def\dist{{\rm dist}}

\newtheorem{corollary}{Corollary}[section]
\newtheorem{theorem}{Theorem}
\newtheorem{problem}{Problem}

\usepackage{booktabs}

\renewcommand{\natural}{{\mathbb{N}}}

\newcommand{\real}{{\mathbb{R}}}

\newcommand{\map}[3]{#1: #2 \rightarrow #3}

\newcommand{\norm}[1]{\|#1\|}

\title{Distributed $n$-player approachability and consensus\\ in 
  coalitional games\thanks{A preliminary conference version of this paper has appeared as  \cite{DB-GN:12Necsys}. The current paper includes, in addition: i) more
  detailed and revised proofs of the main results, ii) analysis of adversarial
  disturbances; iii) analysis of the connections with approachability theory in
  its strategic version for two-player repeated games, and iii) numerical
  studies.  The work of D. Bauso was supported by the 2012 ``Research Fellow'' Program of the Dipartimento di Matematica, Universit\`a di Trento and by PRIN 20103S5RN3 ``Robust decision making in markets and organizations, 2013-2016''. The second author wants to thank J. Hendrickx for the helpful discussion on the
  proof of Theorem~\ref{thm2}.}
} 

\author{Dario Bauso\thanks{D. Bauso is with Dipartimento di Ingegneria Chimica, Gestionale, Informatica e Meccanica,
Universit\`a di Palermo, Italy, email: \textsl{dario.bauso@unipa.it}. D. Bauso is currently academic visitor at the Department of Engineering Science, University of Oxford, UK.} and
Giuseppe Notarstefano\thanks{G. Notarstefano is with Department
    of Engineering, Universit\`a del Salento (University of Lecce), Via per Monteroni, 73100, Lecce, Italy
    \textsl{giuseppe.notarstefano@unisalento.it}}
}

\begin{document}


%

\maketitle

\begin{abstract}                          
We study a distributed allocation process where, repeatedly in time, every player renegotiates past allocations with neighbors and allocates new revenues. The average allocations evolve according to a doubly (over time and space) averaging algorithm. We study conditions under which the average allocations reach consensus to any point within a predefined target set even in the presence of adversarial disturbances. Motivations arise in the context of coalitional games with transferable utilities (TU) where the target set is any set of allocations that make the grand coalitions stable.  

\end{abstract}


\section{Introduction}
\label{sec:introduction}

We consider a two-step distributed allocation process where at every time players first renegotiate their past allocations and second generate a new revenue and allocate it. The time-averaged allocations evolve according to a \emph{doubly (over time and space) averaging dynamics}. The goal is to let all allocations reach consensus to any value in a predefined set even in the presence of an adversarial disturbance. 

\noindent 
{\bf Motivations.} The problem arises in the context of dynamic coalitional games with Transferable Utilities (TU games) \cite{NB13}. A coalitional TU game consists in a set
of players, who can form coalitions, and a characteristic function that provides
a value for each coalition. The predefined set introduced above can be thought of as (but it is not limited to) the core of the game. This is the set of imputations under which no coalition has a value greater than the sum of its members' payoffs. Therefore, no coalition has incentive to leave the grand coalition and receive a larger payoff.

%

\noindent
{\bf Highlights of contributions.} 
We analyze conditions under which  the average allocations: (i)  \emph{approach} the set $X$ (Theorem~\ref{thm1}), (ii)  reach \textit{consensus}, in which case we also compute the consensus value (Theorem~\ref{thm2}), and (iii) are robust against disturbances (Theorem~\ref{thm3}). 
  
\noindent {\bf Related literature. } Coalitional games with transferable utilities (TU) were first introduced by von Neumann and Morgenstern~\cite{VM}. 
Here, a main issue is to study whether the core is
an ``approachable'' set, and which allocation processes can drive the ``complaint vector'' to that set. Approachability theory was developed by Blackwell in the early '56, \cite{B56},
and is captured in the well known Blackwell's Theorem. 
The geometric (approachability) principle that lies behind the Blackwell's
Theorem is among the fundamentals in  allocation processes in coalitional games \cite{L02-CT}.  
The discrete-time dynamics analyzed in the paper follows the rules of a typical
consensus dynamics (see, e.g., \cite{NOP10} and references
therein). 
among multiple agents, where an underlying communication graph for the agents
and balancing weights have been used with some
variations to reach an agreement on common decision
variable in~\cite{Nedic09,NOOT09,NOP10,Sundhar09,GN-FB:07z,MB-GN-FB-FA:12} 
for distributed multi-agent optimization.

The paper is organized as follows. 
In Section, \ref{sec:algorithm}, we formulate the problem and discuss motivations and main assumptions. In Section \ref{results}, we illustrate the main results. In Section \ref{sims} we provide numerical illustrations. Finally, in Section \ref{conc}, we provide concluding
remarks and future directions.

\noindent {\bf Notation}.  We view vectors as columns. 
For a vector $x$, we use
$[x]_j$ to denote its $j$th coordinate component.  
We let $x'$ denote the transpose of a
vector $x$, and $\|x\|$ denote its Euclidean norm.  An $n\times n$ matrix $A$ is
row-stochastic if the matrix has nonnegative entries $a^i_j$ 
and $\sum_{j=1}^n
a^i_j=1$ for all $i=1,\ldots,n$. For a matrix $A$, we use $a^i_j$ or
$[A]_{ij}$ to denote its $ij$th entry. A matrix $A$ is doubly stochastic if both
$A$ and its transpose $A'$ are row-stochastic.  Given two sets $U$ and $S$, we
write $U\subset S$ to denote that $U$ is a proper subset of $S$.  We use $|S|$
for the cardinality of a given finite set $S$.
We write $P_X[x]$ to denote the projection of a vector $x$ on a set $X$, and we
write $\dist(x,X)$ for the distance from $x$ to $X$, i.e., $P_X[x] = \arg
\min_{y \in X} \|x - y\|$ and $\dist(x,X)=\|x-P_X[x]\|$, respectively.  
%
Given a function of time $x(\cdot):\mathbb N \rightarrow \mathbb R$, we denote by
$\bar x(t)$ its average up to time $t$, i.e., $\bar x(t):=\frac{1}{t}
\sum_{\tau=1}^t x(\tau)$.


\section{Distributed reward allocation algorithm}
\label{sec:algorithm}
Every player in a set $N =\{1,\ldots,n\}$ is characterized by an average allocation vector $\hat x_i(t+1)\in \mathbb R^n$. At every time he renegotiates with \emph{neighbors} all past allocations and generates a new allocation vector $x_i(t+1)$.  
The time-averaged allocation $\hat{x}_i (t)$ evolves as follows: 
\begin{equation}\label{dyn}
  \hat x_i(t+1) = \frac{t}{t+1} \left[\sum_{j=1}^n a^i_j(t) \hat x_j(t)\right] +
  \frac{1}{t+1} x_i(t+1),
\end{equation}
where $a^i=(a^i_1,\ldots,a^i_n)'$ is a vector of nonnegative weights consistent
with the sparsity of the \emph{communication graph} $\mathcal{G}(t)=(N,\mathcal{E}(t))$.  A link $(j,i)\in \mathcal{E}(t)$ exists if player $j$ is a neighbor of player~$i$ at time $t$, i.e. if player $i$ renegotiates allocations with player~$j$ at time $t$.




%



{\bf Problem.}
Our goal is to study under what conditions all allocation vectors converge to
a unique value and this value belongs to a predefined set $X$: for all $i,j \in V$,
\begin{equation}
\hat{x}_i(t)=\hat{x}_j(t) \in X,\quad \mbox{for $t\rightarrow
    \infty$}.
\end{equation} 

\noindent
In the sequel, we rewrite equation (\ref{dyn}) in the compact form:
\begin{equation}
\label{dyn1} 
\hat x_i(t+1) = \frac{t}{t+1} w_i(t) + \frac{1}{t+1} x_i(t+1),
\end{equation}
where $w_i(t)$ is the \textit{space average} defined as
\begin{equation}
\label{wik}
w_i(t)=\left[\sum_{j=1}^n a^i_j(t) \hat{x}_j(t)\right].
\end{equation}


\subsection{Motivations}\label{sec:scenario}
The set $X$ introduced above can be thought of  as  the core of a coalitional game with Transferable Utilities (TU game). 

A coalitional TU game is defined by a pair $<N,\eta>$, where $N =
\{1,\ldots,n\}$ is a set of players and $\map{\eta}{2^N}{\real}$ a function
defined for each coalition $S\subseteq N$ ($S\in 2^N$). The function $\eta$
determines the value $\eta(S)$ assigned to each coalition $S\subset N$, with
$\eta(\emptyset)=0$.
We let $\eta_S$ be the value $\eta(S)$ of the characteristic function
$\eta$ associated with a nonempty coalition $S\subseteq N$.  Given a TU game
$<N,\eta>$, let $C(\eta)$ be the core of the game,
\[
\begin{split}
 C(\eta) &=\left\{x
   \in\mathbb{R}^{n} \,\Big|\, \sum_{j\in N}[x]_j=\eta_N,\ \right.\\ 
 &\qquad\left. \sum_{j\in S}
   [x]_j\ge\eta_S\hbox{ for all nonempty } S\subset N\right\}.
\end{split}
\]
Essentially,  the core of the game is the set of all allocations that make the
grand coalition stable with respect to all subcoalitions.
Condition $\sum_{j\in N} [x]_j=\eta_N$ is also called efficiency
condition. Condition $\sum_{j\in S} [x]_j\ge\eta_S$ for all nonempty $S\subset
N$ is referred to as ``stability with respect to subcoalitions'', since it
guarantees that the total amount given to the members of a coalition exceeds the
value of the coalition itself.


\subsection{Main assumptions}
Following 
\cite{NOP10} (see also \cite{NB13}) we can make the following assumptions on the
information structure.
We let $A(t)$ be the weight matrix with entries $a^i_j(t)$.

\begin{assumption}\label{assum:weights}
  Each matrix $A(t)$ is doubly stochastic with positive diagonal.  Furthermore,
  there exists a scalar $\alpha>0$ such that $a^i_j(t)\ge \alpha$ whenever
  $a^i_j(t)>0$.
\end{assumption}


At any time, the instantaneous graph $\mathcal{G}(t)$ need not be
connected. However, for the proper behavior of the process, the union of the
graphs $\mathcal{G}(t)$ over a period of time is assumed to be connected. 

\begin{assumption}\label{assum:graph}
  There exists an integer $Q\ge1$ such that the graph
  $\left(N,\bigcup_{\tau=tQ}^{(t+1)Q-1}\mathcal{E}(\tau)\right)$ is strongly
  connected for every $t\ge0$.
\end{assumption}

It is worth noting that the above assumptions are fairly standard in the
  distributed computation literature. In particular, the joint strong
  connectivity is the weakest possible assumption to guarantee persistent
  circulation of the information through the graph. The double stochasticity of
  the matrix $A(t)$ is a common assumption to guarantee average consensus.


Let $X \subset \mathbb R^n$ be the core set of the game. A common assumption in
approachability theory is that both the core set is convex and bounded, and the
payoff (or loss) vectors generated at each time are bounded. Thus, following
\cite{B56,CL06}, we borrow and adapt such an assumption to our framework.


\begin{assumption}
\label{asm:unit}
The core set $X$ is nonempty.
\end{assumption}

Notice that a nonempty core is a convex and compact set.


The next assumption indicates how the new reward vector has to be generated in order
to obtain approachability.
\begin{assumption}
\label{asm:app}
For each $i\in N$ the new reward vector $x_i(\cdot)$ is bounded, i.e., there exists
$L>0$ s.t. $\forall t\geq0 \; \|x_i(t+1)\| \leq L$, and satisfies the following inequality, for a scalar
  negative number, $\phi <0$,
\[
\left(w_i(t)-P_X(w_i(t))\right)' \left(x_i(t+1) - P_X(w_i(t)\right)\leq \phi <
0.
\]
\end{assumption}
From a geometric standpoint, Assumption \ref{asm:app} requires that, given the
two half-spaces identified by the supporting hyperplane of $X$ through
$P_X(w_i(t))$, the new reward vector $x_i(t+1)$ lies in the half-space not containing
$w_i(t)$.

 \begin{figure} 
   \centering
 \includegraphics[width=0.5\linewidth]{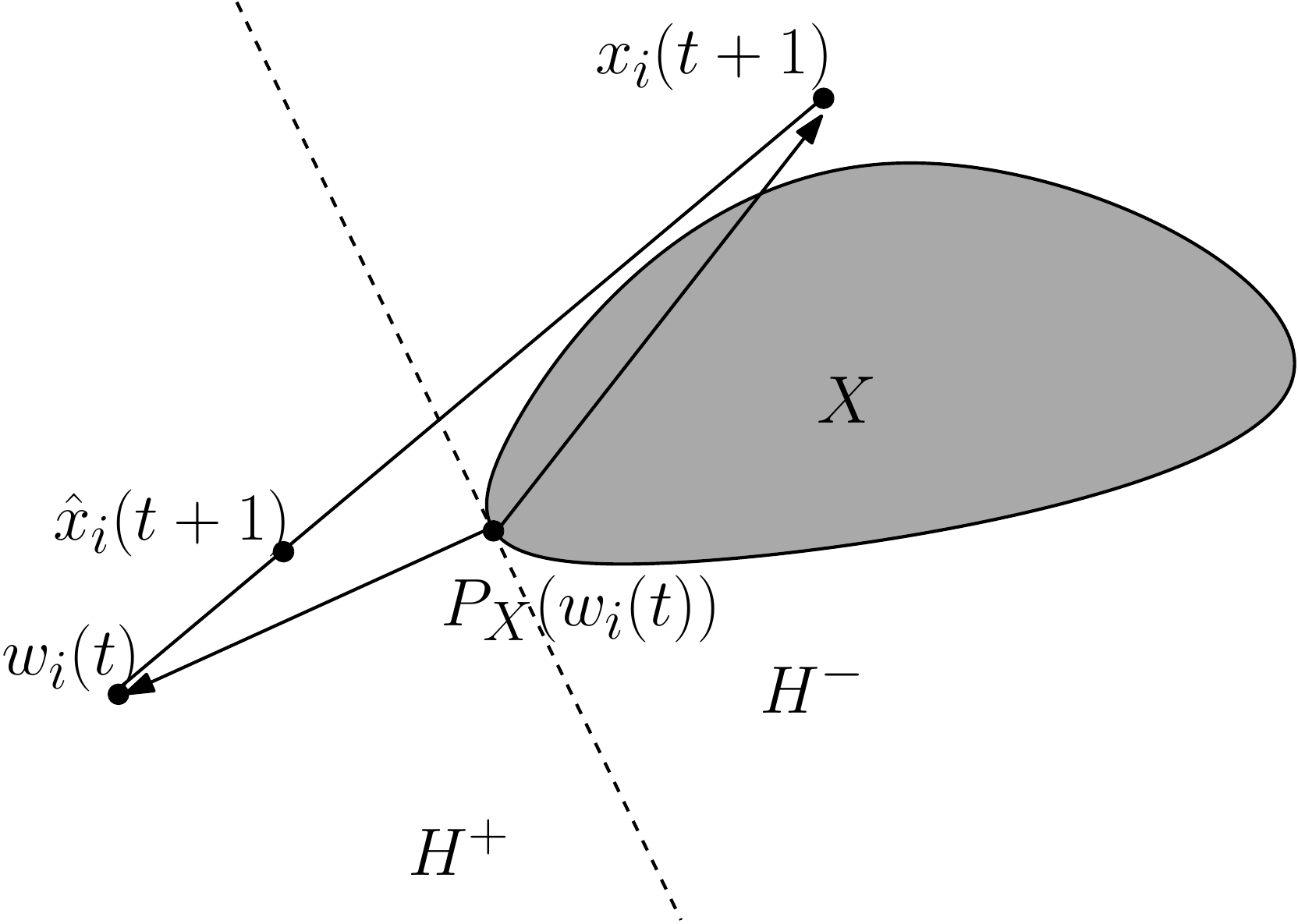}
   \caption{Approachability principle.}\label{fig:graphs}
 \end{figure}





\section{Main results}
\label{results}
Next, we provide the main results of the paper. Namely, we prove that  the
average allocations: (i)  \emph{approach} the set $X$ (Theorem~\ref{thm1}), (ii)  reach \textit{consensus}  (Theorem~\ref{thm2}), and (iii) are robust against disturbances (Theorem~\ref{thm3}).
\subsection{Approachability and consensus}
Before stating the first theorem, we need to introduce two lemmas. The next lemma
establishes that the space averaging step in (\ref{dyn}) reduces the total
distance (i.e. the sum of distances) of the estimates from the set $X$.
\begin{lemma}\label{lem1}
  Let Assumption~\ref{assum:weights} hold. Then the total distance from $X$
  decreases when replacing the allocations $\hat x_i(t)$ by their space averages
  $w_i(t)$, i.e.,
\[
  \sum_{i=1}^{n} \dist(w_i(t),X) \leq \sum_{i=1}^{n} \dist(\hat x_i(t),X).
\]
\end{lemma}
%
%
As a preliminary step to the next result, observe that, from the definition of
$\dist(\cdot,X)$ and from (\ref{dyn}) and (\ref{wik}), it holds 
\begin{equation}
\begin{split}
\label{eq:lem2-1}
   \dist&(\hat x_i(t+1),X)^2  = \|\hat x_i(t+1)-P_X[\hat x_i(t+1)]\|^2 \\ 
    &\leq  \|\hat x_i(t+1)-P_X[w_i(t)]\|^2\\ 
    &=  \left \|\frac{t}{t+1} w_i(t) + \frac{1}{t+1} x_i(t+1)-P_X[w_i(t)] \right \|^2 \\ 
    &=  \left \|\frac{t}{t+1} \left(w_i(t) - P_X[w_i(t)]  \right)\right.\\    
     & \qquad\left.+  \frac{1}{t+1} \left(x_i(t+1) - P_X[w_i(t)]  \right)  \right \|^2\\
    &=  \left(\frac{t}{t+1}\right)^2 \|w_i(t) - P_X[w_i(t)]\|^2\\
   & \qquad  + \left( \frac{1}{t+1}\right)^2 \|x_i(t+1) - P_X[w_i(t)] \|^2 \\ 
 & +  \frac{2t}{(t+1)^2} (w_i(t) - P_X[w_i(t)])' (x_i(t+1)-P_X[w_i(t)]).
\end{split}
\end{equation}

The following lemma states that, under the approachability assumption, the
distance of each single estimate from $X$ decreases with respect to the one of
the spatial average when applying the time averaging step.
\begin{lemma}\label{lem2}
  Let Assumptions \ref{asm:unit}-\ref{asm:app} hold. Then, there exists a
  positive integer scalar, $\tilde{t} >0$, such that for all $t \geq \tilde t
  >0$ the distance of each single $\hat x_i(t+1)$ decreases in comparison with
  the distance of $w_i(t)$, i.e.,
\[
\dist(\hat x_i(t+1),X) < \dist(w_i(t),X),\quad \forall i=1,\ldots,n
\]
\end{lemma}

We are now ready to state the first main result.
\begin{theorem}
\label{thm1} 
Let Assumptions \ref{assum:weights}-\ref{asm:app} hold.  Then all average allocations
approach set $X$, i.e.,
\[
  \lim_{t \rightarrow \infty} \sum_{i=1}^{n} \dist(\hat x_i(t),X) = 0.
\]
\end{theorem}

Next, let us introduce the barycenter of respectively the estimates and the reward vectors 
\[
\hat{x}_b(t) :=\frac{1}{n}
\sum_{i=1}^n \hat{x}_i(t)
\qquad \text{and} \qquad
x_b(t):=\frac{1}{n}
\sum_{i=1}^n x_i(t).
\] 
%
Consistently, let us denote as $\bar{x}_b(t)$ the time average of the
barycenter, i.e. 
\[
\bar{x}_b(t) = \frac{1}{t+1} \sum_{\tau=0}^t x_b(\tau).
\]
The following lemma establishes that the barycenter of the estimates evolves as
the time average $\bar{x}_b(t)$ of the barycenter of the reward vectors
generated by the players.

\begin{lemma}\label{lem3}
  The barycenter of the local allocations $\hat{x}_b(t)$ coincides at each time $t$ with the
  time-average of the barycenter of the generated reward vectors $\bar{x}_b(t)$.
\end{lemma}
%
The following theorem establishes that all allocations converge to $\bar{x}_b(t)$,
which in the limit must belong to $X$ according to Theorem~\ref{thm1}.
\begin{theorem}\label{thm2}(Consensus to the barycenter time-average)
  Let Assumptions \ref{assum:weights}-\ref{asm:app} hold.  Then, all players
  reach consensus on the time-average of the barycenter of the reward vectors
  generated by each player, $\bar{x}_b(t)$, i.e.,
\[
\lim_{t \rightarrow \infty} \|\hat x_i(t)-\bar{x}_b(t)\|=0 \quad \forall i=1,\ldots,n.
\]
\end{theorem}

Summarizing the two main results, we have proven that asymptotically all the
players' allocations converge to the time-average of the barycenter of the
generated reward vectors and that this vector lies in the core of the game.

\subsection{Adversarial disturbance}\label{zsg}
Here we analyze the case where, for each player $i\in N$, the input $x_i(\cdot)$ is the payoff of a
repeated two-player game between player $i$  (Player $i_1$) and an (external) adversary
(Player $i_2$). With some slight abuse of notation we denote $S_1$ and $S_2$ the
finite set of actions of players $i_1$ and $i_2$ respectively.

The instantaneous payoff $x_i(t)$ at time $t$ is given by a function
$\phi_i: S_1 \times S_2 \to \mathbb R^n$ as follows:
\[
 x_i(t)= \phi(j(t),k(t)),
\]
where  $j(t) \in S_1$ and $k(t) \in S_2$.
We extend $x_i$ to the set of mixed actions pairs, $\Delta(S_1) \times
\Delta(S_2)$, in a bilinear fashion.  In particular, for every pair of mixed
strategies $(p(t),q(t)) \in \Delta(S_1) \times \Delta(S_2)$ for player $i_1$ and $i_2$
at time $t$, the expected payoff is 
\[
\mathbb E x_i(t)=\sum_{j \in S_1} \sum_{k
  \in S_2} p_j(t) q_k(t) \phi(j,k).
\]
For simplicity the one-shot vector-payoff
game $(S_1,S_2,x_i)$ is denoted by $G_i$.

Let $\lambda\in \mathbb{R}^n$. Denote %
 by $\langle \lambda, G_i \rangle$ the zero-sum one-shot game whose set of players and their
action sets are as in the game $G_i$, and the payoff
that player 2 pays to player 1 is $\lambda' \phi(j,k)$ for every  $(j,k)\in S_1 \times S_2$.

The resulting zero-sum game is described by the matrix 
\[
  \Phi_\lambda=[\lambda' \phi(j,k)]_{j\in S_1,k\in S_2}.
\]
As a zero-sum one-shot game, the game  $\langle \lambda, G_i \rangle$ has a value, denoted 
\[
v_{\lambda}:=\min_{p \in \Delta S_1} \max_{q \in \Delta S_2} p' \Phi_\lambda q =
\max_{q \in \Delta S_2} \min_{p \in \Delta S_1} p' \Phi_\lambda q.
\]

For every mixed action $p \in \Delta(S_1)$ denote $D_1(p)$ the set of all payoffs that might be realized when
player $i_1$ plays the mixed action $p$:
\[ D_1(p) = \{ x_i(p,q) \colon q \in \Delta(S_2)\}. \]
If $v_{\lambda}\geq 0 $ (resp.
$v_{\lambda}>0 $), then there is a mixed action $p \in \Delta(S_1)$
such that $D_1(p)$ is a subset of the closed half space $\{x\in
\mathbb{R}^n\colon  \lambda' x \geq 0\}$ (resp. half
space $\{x\in \mathbb{R}^m\colon \lambda' x  >0\}$).


Let us introduce next the counterpart of Assumption \ref{asm:app} in this new worst-case setting.
\begin{assumption} 
\label{asm:app_wc}
For any $w_i(t) \in \mathbb R^n$, there exists a mixed strategy $p(t+1) \in
\Delta(S_1)$ for Player $i_1$ such that, for all mixed strategy $q(t+1)\in
\Delta(S_2)$ of Player $i_2$, the new reward vector $x_i(\cdot)$ is bounded,
i.e. there exists $L>0$ s.t. $\forall t\geq0 \; \|x_i(t+1)\| \leq L$, and
satisfies
\[
\left(w_i(t)-P_X(w_i(t))\right)' \left(\mathbb E x_i(t+1) - P_X(w_i(t)\right)\leq \phi <
0,
\]
where $\mathbb E  x_i(t+1)=\sum_{j \in S_1} \sum_{k \in S_2} p_j(t+1) q_k(t+1) \phi(j,k).$
\end{assumption}
The above condition is among the foundations of approachability theory as it
guarantees that the average payoff $\frac{1}{T} \sum_{t=0}^{T-1} x_i(t)$
converges almost surely to $X$ (see, e.g., \cite{B56} and also \cite{CL06},
chapter 7). Here we adapt the above condition to the multi-agent and distributed
scenario under study.

\begin{corollary}[see \cite{B56}, Corollary 2]
Any convex set $X \subset \mathbb R^n$ is approachable if and only if $v_\lambda < 0$ for any $\lambda \in \mathbb R^n$.
\end{corollary}
Next we show that if the approachability condition expressed above holds true,
then $dist(\hat x_i,X)$ tends to zero for any $X$.  We write $w.p.1$ to mean
``with probability 1''.

\begin{theorem}
\label{thm3} 
Let Assumptions \ref{assum:weights}-\ref{asm:unit} and \ref{asm:app_wc} hold.  Then all average allocations
approach set $X$, i.e.,
\[
  \lim_{t \rightarrow \infty} \sum_{i=1}^{n} \dist(\hat x_i(t),X) = 0, \quad w.p.1.
\]
\end{theorem}

We conclude this section by observing that Theorem \ref{thm2} still holds and
therefore all players' estimates reach consensus on the time-average of the
barycenter of the reward vectors generated by each player.

\section{Simulations}\label{sims}
We illustrate the results in a game with four players, $N=\{1,\ldots, 4\}$, communicating according to a fixed
undirected cycle graph. That is, $\mathcal{G}(t)=(N,\mathcal{E})$ where
$\mathcal{E} = \{(i,j) \;|\; j=i+1, i\in \{1,\ldots n-1\} \; \text{or} \; (i,j)
= (n,1)\}$. 

%

We set $\eta_{\{1\}}=\ldots=\eta_{\{4\}}=2$, $\eta_{\{1,2\}} = 5$, $\eta_{\{3,4\}} =
5$, $\eta_{\{1,2,3\}} = 7$ and $\eta_N = 10$ ($\eta_S$ is the value of coalition $S$). That is, each player expects to receive at least a reward of $2$ which is its value as a singleton
coalition. But, for example, players $1$ and $2$ expect to be more valuable if
they form a coalition as well as $3$ and $4$. 
Consistently, the core of the game is the polyhedral set given by
\[
\begin{split}
C(\eta) =\Big\{x
  \in\mathbb{R}^{4} \,\Big|\, x_1+x_2+x_3+x_4=10,\\ x_1+x_2+x_3 \geq 7, \
  x_1+x_2\geq 5,\\
\ x_3+x_4\geq 5,
x_1\geq 2, \ldots, x_4\geq 2\Big\}.
\end{split}
\] 

\begin{figure}[h]
  \centering
  \includegraphics[width=0.48\linewidth]{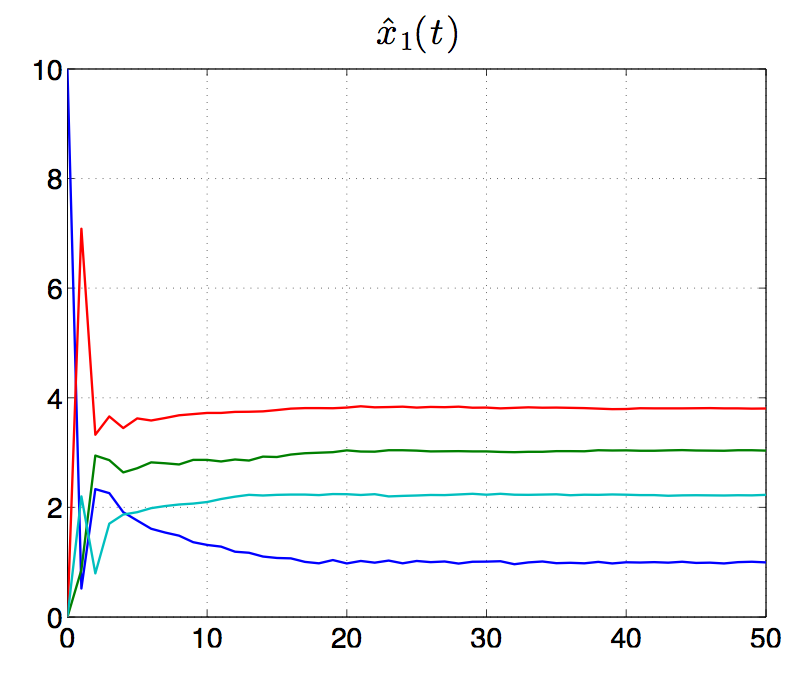}\,
  \includegraphics[width=0.48\linewidth]{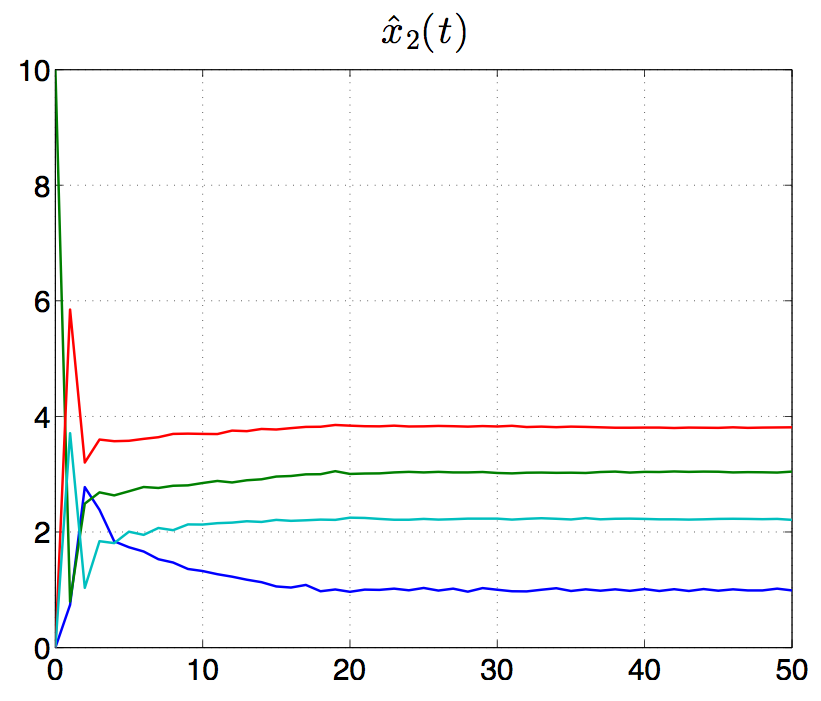}\\[1.2ex]
  \includegraphics[width=0.48\linewidth]{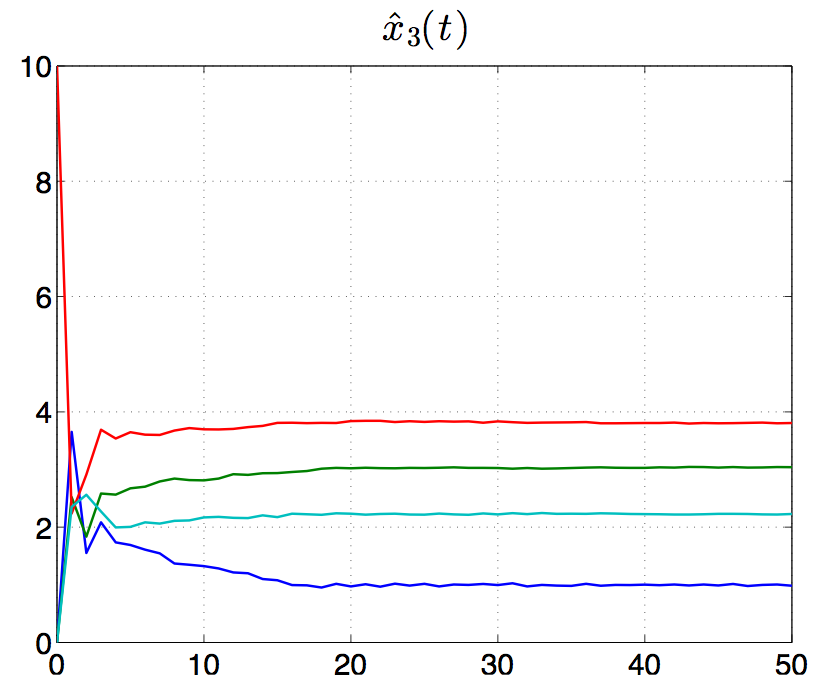}\,
  \includegraphics[width=0.48\linewidth]{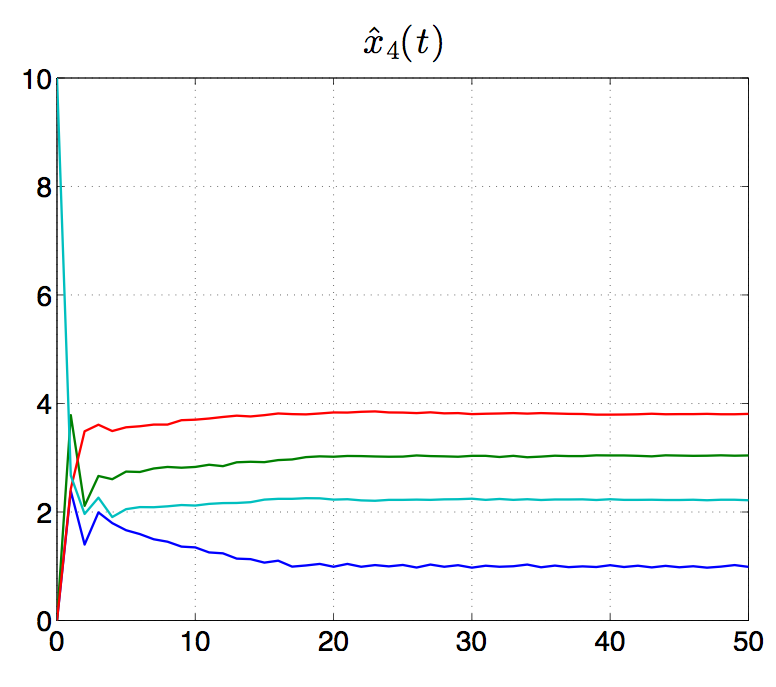}
  \caption{Local average reward vectors}
  \label{fig:rewards}
\end{figure}

We initialize the assignments assuming each player assign itself the entire
reward.  That is, denoting $b_i \in \real^n$ the $i$-th canonical vector (so
that, e.g., $b_1 = [1 \; 0 \; \ldots \; 0]'$), we set $\hat{x}_i(0) = 10 \,b_i$
for all $i\in\{1,\ldots,n\}$. At every iteration $t\in \natural$, each player
chooses the new reward vector $x_i(t+1)$ according to the approachability
principle. In particular, we set $x_i(t+1) = P_X[w_i(t)] +
\alpha\, (P_X[w_i(t)]-w_i(t)) + v^{\top}$, where $\alpha$ is a random number
uniformly distributed in $[0,1]$ and $v^{\top}$ a random vector belonging to the
hyperplane tangent to the core at $P_X[w_i(t)]$ with coordinates uniformly
chosen in $[0,1]$. The temporal evolution of the local estimates of the average
reward vector is depicted in Figure~\ref{fig:rewards}. As expected the local
estimates converge to the same average assignment which is the point of the core
$[ 3.8 \;\, 3 \;\, 2.2 \;\, 1]'$.

\section{Conclusions}\label{conc}

We have analyzed convergence conditions of a distributed allocation process arising in the context of TU games. Future directions include the extension of our results  to population games with mean-field interactions, and averaging algorithms driven by Brownian motions.

%

\bibliographystyle{plain} 
\bibliography{biblio}

\section*{Appendix}
\subsection*{Proof of Lemma \ref{lem1}}
By convexity of the distance function $\dist(\cdot,X)$ and from (\ref{wik}) we have
\[
\dist(w_i(t),X) \leq \sum_{j=1}^n a^i_j(t) \dist(\hat x_j(t),X).
\]
Summing over $i=1,\ldots,n$ both sides of the above inequality we obtain
\[
\begin{split}
  &\sum_{i=1}^n \dist(w_i(t),X) \leq \sum_{i=1}^n \sum_{j=1}^n a^i_j(t) \dist(\hat
  x_j(t),X) \\= &\sum_{j=1}^n \left(\sum_{i=1}^n a^i_j(t)\right) \dist(\hat x_j(t),X)=
  \sum_{j=1}^n \dist(\hat x_j(t),X),\\
\end{split}
\] 
where the last equality follows from the stochasticity of $A(t)$ in
Assumption~\ref{assum:weights}. This concludes the proof.

\subsection*{Proof of Lemma \ref{lem2}}
Rearranging equation (\ref{eq:lem2-1}) we obtain
\begin{equation}
\begin{split}
\label{eq1} 
\|\hat x_i(t+1)-P_X[\hat x_i(t+1)]\|^2 \\
- \frac{t^2}{(t+1)^2} \|w_i(t)-P_X[w_i(t)]\|^2 \leq\\
\; \frac{1}{(t+1)^2}\|x_i(t+1) - P_X[w_i(t)] \|^2 \\
+  \frac{2t}{(t+1)^2} (w_i(t) - P_X[w_i(t)])'
  (x_i(t+1)-P_X[w_i(t)]).\\
\end{split}
\end{equation}
Note that the left hand side in (\ref{eq1}) approximates $\dist(\hat x_i(t+1),X)^2 -
\dist(w_i(t),X)^2$ for increasing $t$ and also that for all $t$ the left hand side upper
bounds such a difference, i.e., 
\begin{equation*}
\begin{split}
\dist(\hat x_i(t+1),X)^2 - \dist(w_i(t),X)^2 \\ \leq
 \dist(\hat x_i(t+1),X)^2 - \frac{t^2}{(t+1)^2} \dist(w_i(t),X)^2  \quad \forall t.
\end{split}
\end{equation*}

It remains to note that there exists a great enough scalar integer $\tilde
t$ 
such that the left hand side in (\ref{eq1}) is negative for all $t \geq
\tilde t$.
From the boundedness of set $X$ and of vectors $x_i(t)$, there exists
$M>0$ such that $\|x_i(t+1) - P_X[w_i(t)] \|^2<M$. Thus, we have
\begin{equation}
\begin{split}
\dist(\hat x_i(t+1),X)^2 - \dist(w_i(t),X)^2 \\ \leq 
 \dist(\hat x_i(t+1),X)^2 - \frac{t^2}{(t+1)^2} \dist(w_i(t),X)^2 \\
  \leq \frac{1}{(t+1)^2}\big(\|x_i(t+1) - P_X[w_i(t)] \|^2 \\
+ 2t (w_i(t) - P_X[w_i(t)])' (x_i(t+1)-P_X[w_i(t)])\big)\\
  \leq \frac{1}{(t+1)^2}(M + 2 t \phi) <0
\end{split}
\end{equation}
Taking $\tilde{t}> -M/2\phi >0$ concludes the proof. 

\subsection*{Proof of Theorem \ref{thm1}}
Recall from (\ref{eq:lem2-1}) that
\begin{equation*}
\begin{split}
  \|\hat x_i(t+1)-P_X[\hat x_i(t+1)]\|^2  \leq  \\
   \left(\frac{t}{t+1}\right)^2 \|w_i(t) - P_X[w_i(t)]\|^2\\
    + \left( \frac{1}{t+1}\right)^2 \|x_i(t+1) - P_X[w_i(t)] \|^2 \\ 
   \;+  2 \frac{t}{(t+1)^2} (w_i(t) - P_X[w_i(t)])' (x_i(t+1)-P_X[w_i(t)]).
\end{split}
\end{equation*}
From Lemma \ref{lem1} and rearranging the above inequality, we have
\begin{equation*}
\begin{split}
  \sum_{i=1}^n \big[(t+1)^2 \|\hat x_i(t+1)-P_X[\hat
    x_i(t+1)]\|^2 \\- t^2 \|\hat x_i(t)-P_X[\hat x_i(t)]\|^2\big] \\[1.2ex]
   \leq\sum_{i=1}^n \Big[\|x_i(t+1) - P_X[w_i(t)] \|^2 \\+ 2 t (w_i(t) -
    P_X[w_i(t)])' (x_i(t+1)-P_X[w_i(t)])\Big]\\[1.2ex]
  \leq \sum_{i=1}^n \left[\|x_i(t+1) - P_X[w_i(t)] \|^2 \right. ,
\end{split}
\end{equation*}
where the last inequality is due to Assumption~\ref{asm:app}.  Summing
over $t=0,\ldots,\tau-1$, and noting that $\|x_i(t+1) - P_X[w_i(t)]\|$ is
bounded (from Assumption \ref{asm:unit}), so that the right hand side is upper
bounded by some $M>0$, we obtain
\[
\sum_{i=1}^n \tau^2 \|\hat x_i(\tau)-P_X[\hat x_i(\tau)]\|^2 \leq M \tau
\]
from which $\|\hat x_i(\tau)-P_X[\hat x_i(\tau)]\|^2 \leq \frac{M}{\tau}$, and therefore
$\lim_{\tau\rightarrow \infty} \|\hat x_i(\tau)-P_X[\hat x_i(\tau)]\|^2=0$,
which concludes the proof.

\subsection*{Proof of Lemma \ref{lem3}}
  To prove the statement observe that $\bar{x}_b(0) = \hat{x}_b(0) =
  x_b(0)$. Thus, we prove that $\bar{x}_b(t)$ and $\hat{x}_b(t)$ satisfy the
  same dynamics.
  By definition of time-average, $\bar{x}_b(t)$ satisfies the dynamics
  \begin{equation}
    \label{eq:aver_baryc}
    \bar{x}_b(t+1) = \frac{t}{t+1} \bar{x}_b(t) + \frac{1}{t+1} x_b(t+1).
  \end{equation}
  The dynamics of $\hat{x}_b(t)$ is
 \begin{eqnarray*}
 \frac{1}{n} \sum_{i=1}^n\hat{x}_i(t+1) = \frac{1}{n} \Big[ \frac{t}{t+1} \sum_{i=1}^n \sum_{j=1}^n a^i_j(t)
  \hat{x}_j(t) \\+ \frac{1}{t+1} \sum_{i=1}^n x_i(t+1)\Big].
\end{eqnarray*}
  Exchanging the sum signs
  \[
  \hat{x}_b(t+1) = \frac{1}{n}  \frac{t}{t+1} \sum_{j=1}^n \sum_{i=1}^n a^i_j(t) \hat{x}_j(t)
  + \frac{1}{t+1} x_b(t+1),
  \]
  and, by Assumption~\ref{assum:weights} ($A(t)$ is doubly stochastic),
 \begin{equation*}
\begin{split}
  \hat{x}_b(t+1)  =  \frac{1}{n} \frac{t}{t+1} \sum_{j=1}^n \hat{x}_j(t) + \frac{1}{t+1}
  x_b(t+1)\\
  =  \frac{t}{t+1}  \hat{x}_b(t) + \frac{1}{t+1} x_b(t+1),
  \end{split}
  \end{equation*}
  which is the same dynamics as \eqref{eq:aver_baryc}, thus concluding the
  proof. 

 \subsection*{Proof of Theorem \ref{thm2}}
  Using the previous lemma we can show that $\hat x_i(t)$ converges to $\hat
  x_b(t)$. Let us introduce the error of the estimate $\hat{x}_i(t)$ from the
  barycenter, i.e. $\hat{e}_i(t) = \hat{x}_i(t) - \hat{x}_b(t)$. The error
  dynamics is given by
  \[
\begin{split}
  \hat{e}_i(t+1) = \frac{t}{t+1} \left[ \sum_{j=1}^n a^i_j(t) \hat{e}_j(t) +
  \sum_{j=1}^n a_j^i \hat{x}_b(t)\right]
   \\ + \frac{1}{t+1} e_i(t+1) + \frac{1}{t+1}
  x_b(t+1)\\[1.2ex] 
 - \frac{t}{t+1} \hat{x}_b(t) - \frac{1}{t+1}x_b(t+1),
\end{split}  
\]
where $e_i(t) = x_i(t) - x_b(t)$. Thus
  \[
  \hat{e}_i(t+1) = \frac{t}{t+1} \big( \sum_{j=1}^n a^i_j(t) \hat{e}_j(t)\big) +
  \frac{1}{t+1} e_i(t+1).
  \]
  Multiplying both sides by $(t+1)$ and taking $t$ inside the sum,
  \[
  (t+1) \hat{e}_i(t+1) = \sum_{j=1}^n a^i_j(t) t \hat{e}_j(t) + e_i(t+1).
  \]
  Defining $\hat{z}_i(t) = t \,\hat{e}_i(t)$, we have
  \[
  \hat{z}_i(t+1) = \sum_{j=1}^n a^i_j(t) \hat{z}_j(t) + e_i(t+1).
  \]
In vector form the above equation turns to be
\begin{equation}
  \label{eq:z_dynamics}
    \hat{z}(t+1) = \big(A(t) \otimes I_n \big)\hat{z}(t) + e(t+1),
  \end{equation}
  with $\hat{z}(t) = [z_1(t) \; \ldots \; z_n(t)]'$, $\hat{e}(t) = [e_1(t)
    \; \ldots \; e_n(t)]'$, $I_n$ the identity matrix of dimension $n$ and
    $\otimes$ the Kronecker product. Notice that denoting $[\hat{z}]_\ell = \big[
    [\hat{z}_1]_\ell \; \ldots \; [\hat{z}_n]_\ell\big]$ and $[e]_\ell = \big[[e_1]_\ell
    \; \ldots \; [e_n]_\ell\big]$, $\ell \in \{1,\ldots, n\}$, the dynamics of each
    $[\hat{z}]_\ell$ is given by
    \begin{equation}
      \label{eq:z_dynamics_scalar}
      [\hat{z}]_\ell(t+1) = A(t) [\hat{z}]_\ell(t) + [e]_\ell(t+1).
  \end{equation}
Thus, we can simply work on each component separately. Slightly abusing notation we
neglect the subscript of $[\hat{z}]_\ell$ and $[e]_\ell$, and write $\hat{z}(t)$
and $e(t)$.
  
It is worth noting that the driven system \eqref{eq:z_dynamics_scalar}, and so
\eqref{eq:z_dynamics}, is \emph{not} bounded-input-bounded-state stable (when a
general input signal is allowed). That is, for general initial condition and
input signal the state trajectory may diverge. We show that for the special
initial condition ($\hat{z}(t) = 0$ by construction) and class of input signals
($\mathbf{1}' e(t+1) = 0$ by definition) under consideration, the state
trajectories of \eqref{eq:z_dynamics} are bounded.

First, let us observe that, multiplying both sides of \eqref{eq:z_dynamics} by
the vector $\mathbf{1}' = [1\; \ldots \; 1]$, we get
\begin{equation}
  \label{eq:orthog_to_1}  
  \begin{split}    
  \mathbf{1}'\hat{z}(t+1) &= \mathbf{1}' A(t) \hat{z}(t) + \mathbf{1}'
  e(t+1)\\
 &= \mathbf{1}' \hat{z}(t).
  \end{split}
\end{equation}
Since $\hat{z}(0)=0$ by construction, it holds $\mathbf{1}'\hat{z}(t) = 0$ for
all $t\in\natural$. That is, $\hat{z}(t)$ is orthogonal to the vector
$\mathbf{1}$ for all $t$.

Next, we show that the trajectory $\hat{z}(\cdot)$ is bounded. Following
  \cite{blondel2005convergence}, let $P\in \real^{(n-1) \times n}$ be a matrix
defining an orthogonal projection onto the space orthogonal to
span\{$\mathbf{1}$\}. It holds that $P \mathbf{1}=0$ and $\norm{Px}_2 =
\norm{x}_2$ if $x'\mathbf{1}=0$. Thus, from equation \eqref{eq:orthog_to_1} we
have that $\norm{P\hat{z}(t)}_2 = \norm{\hat{z}(t)}_2$ for all $t$. Therefore,
proving boundedness of $\hat{z}(\cdot)$ is equivalent to showing that
$P\hat{z}(\cdot)$ is bounded.
For a given $P$, associated to any $A(t)$ satisfying
  Assumption~\ref{assum:weights}, there exists $\bar{A}(t)$ satisfying $P A(t)
= \bar{A}(t) P$. The spectrum of $\bar{A}(t)$ is the spectrum of $A(t)$ after
removing the eigenvalue $1$.  Multiplying both sides of equation
\eqref{eq:z_dynamics} by $P$, we get
\begin{equation}
  \label{eq:z_dynamics_orthog}
  \begin{split}    
    P\hat{z}(t+1) &= P A(t) \hat{z}(t) + Pe(t+1)\\
      &= \bar{A}(t) P \hat{z}(t) + P e(t+1).
  \end{split}
  \end{equation}
%
  Under Assumptions~\ref{assum:weights} and \ref{assum:graph}, the undriven
    dynamics $y(t+1) = \bar{A}(t) y(t)$ is uniformly exponentially stable, i.e.,
    $||y(t)||<C \rho^t ||y(0)||$ with $C$ and $\rho<1$ independent of
    $y(0)$ and depending only on $n$, $Q$ and $\alpha$ (see Theorem~9.2
    and Corollary 9.1 in \cite{hendrickx2008graphs}). Thus, the state
    trajectories of \eqref{eq:z_dynamics_orthog} are bounded for any bounded
    signal $P e(t+1)$ with $\mathbf{1}'e(t) = 0$.
Since $\mathbf{1}'e(t) = 0$ for all $t$, we have $\norm{Pe(t)}_2=
\norm{e(t)}_2$ for all $t$, which is bounded.
The proof follows by recalling that $\norm{P \hat{z}(t)}_2 =
\norm{\hat{z}(t)}_2$ and that $\hat{z}(t) = t \hat{e}(t)$. 

 \subsection*{Proof of Theorem \ref{thm3}}
From (\ref{eq:lem2-1}), invoking Lemma \ref{lem1} and using Assumption \ref{asm:app_wc}  we have 
\begin{equation*}
\begin{split}
  \sum_{i=1}^n \big[(t+1)^2 \|\hat x_i(t+1)-P_X[\hat
    x_i(t+1)]\|^2 \\- t^2 \|\hat x_i(t)-P_X[\hat x_i(t)]\|^2\big] \\[1.2ex]
    \leq\sum_{i=1}^n \Big[\|x_i(t+1) - P_X[w_i(t)] \|^2 \\+ 2 t (w_i(t) -
    P_X[w_i(t)])' (x_i(t+1)- \mathbb E x_i(t+1))\Big],
\end{split}
\end{equation*}
Summing over $t=0,\ldots,\tau-1$, and noting that $\|x_i(t+1) - P_X[w_i(t)]
\|$ is upper bounded (from Assumption \ref{asm:unit}) by some $M>0$, we obtain
\[
\begin{split}
\sum_{i=1}^n \|\hat x_i(\tau)-P_X[\hat x_i(\tau)]\|^2 \\ \leq \frac{M}{\tau} +
\frac{1}{\tau} \sum_{t=0}^{\tau-1} \sum_{i=1}^n K^i_{t} 
\|x_i(t+1)-\mathbb E x_i(t+1)\|
\end{split}
\]
where $K^i_{t} = \frac{1}{\tau} 2 t \|w_i(t) - P_X[w_i(t)]\|$.
Now, using $\|x_i(t+1)\| \leq L \; \forall t\geq0$ from Assumption
\ref{asm:app_wc} and from (\ref{dyn1}) and (\ref{wik}) we have that $w_i(t)$ is
bounded which in turn implies that $\|w_i(t) - P_X[w_i(t)]\|$ is bounded.  Then,
the second term in the right-hand side is an average of bounded zero-mean
martingale differences, and therefore the Hoeffding-Azuma inequality (together
with the Borel-Cantelli lemma) immediately implies that
$$\lim_{\tau \rightarrow \infty} \sum_{i=0}^n \|\hat x_i(\tau)-P_X[\hat x_i(\tau)]\|^2=0$$ 
which concludes the proof.

\end{document}